\documentclass{article}
\usepackage{amssymb,amsmath,amsthm}

\newtheorem{theorem}{Theorem}

\DeclareMathOperator{\Gal}	{Gal}
\DeclareMathOperator{\cd}	{cd}
\DeclareMathOperator{\Irr}	{Irr}

\newcommand{\ind}[2] {{|{#1}:{#2}|}}
\newcommand{\iso}{\cong}

\newcommand{\FF}{\mathbb F}

\title{Using a Galois connection to compute character degrees}
\author{Mark L.~Lewis and
	   John K.~McVey\\
	\small\sl  Department of Mathematical Sciences, Kent State University\\
	\small\sl Kent, Ohio~~44242\\
	\small\tt lewis@math.kent.edu~~~jmcvey@math.kent.edu}

\begin{document}
\maketitle
\begin{abstract}
Given a Mersenne prime $q$ and a positive even integer $e$, let $F$ and $E$ be the fields of orders $q$ and $q^e$ respectively.  Let $C$ be a cyclic subgroup of $E^\times$ whose index in $E^\times$ is divisible only by primes dividing $q - 1$.  We compute the character degrees of the group $C \rtimes \Gal (E/F)$ by using the Galois connection between the subfields of $E$ and the Galois group $\Gal (E/F)$.
\end{abstract}

\section{Introduction}
In the papers~\cite{mcvey} and~\cite{mcvey2}, the second author worked on determining the character degree sets for the groups of the form $G = C \rtimes \Gal(E/F)$ where $F < E$ were finite fields and ${C < E^\times}$ was any subgroup whose index $\ind {E^\times}C$ is divisible only by the primes that divide $|F^\times|$.  In Theorem~5 of \cite{mcvey2}, it was proved that $\cd (G)$ is the full list of divisors of $\ind EF$, except possibly when three conditions are simultaneously satisfied: (i) $|F|$ is a Mersenne prime, (ii) the degree $\ind EF$ is even, and (iii) $4$ does not divide $|C|$.  However, the character degree set for $G$ when these conditions were mutually satisfied remained unspecified.  The intent of this note is to fill that void.  In particular, we show the corresponding character degree set is the full list of divisors of $\ind EF$, but with $2$ removed.
\medskip

{\noindent\bf Main Theorem. }{\it Fix a Mersenne prime $q$ and an even integer $e > 1$.  Label by $F$ the field $\FF_q$ and by $E$ the field $\FF_{q^e}$.  Take $\pi$ to be the set of primes dividing $q-1$.  Let $\varGamma = \Gal(E/F)$, and fix the subgroup $C\leq E^\times$ under the assumptions that  $\ind {E^\times\!}C$ is a $\pi$-number and $4$ does not divide $|C|$.  If $G = C \rtimes \varGamma$, then
$$\cd (G) = \{n\mid n~\text{divides}~e\}\setminus\{2\}.$$
}

\bigskip\goodbreak

\section{The Degree Computations}

Throughout the following, all groups are assumed to be finite, and standard notations from \cite{Isaacs} are used. Theorem~4 from~\cite{mcvey2} will help in our computations when calculating the stabilizers of irreducible characters.  We quote here the version of it that will be applicable to our situation.  Also, some notation will prove useful.  Given a subset $X \subseteq E$, we denote by $\widehat X$ the smallest subfield of $E$ which contains $X\cup F$.  In particular, a Galois automorphism for $E$ over $F$ centralizes $X$ if and only if it centralizes $\widehat X$.


\begin{theorem}[Theorem 4 of \cite{mcvey2}]\label{MersenneException}
Let $q$ be a Mersenne prime, $e>1$ an even integer, and $\pi$ the set of primes dividing $q-1$.  Fix the fields $F=\FF_q$ and $E=\FF_{q^e}$ and a subgroup $C$ of $E^\times$ whose index $\ind{E^\times\!}C$ is a $\pi$-number.  Assume further that $4$ does not divide~$|C|$.  Then for all fields $F\leq L \leq E$, the equality $L = \widehat{L \cap C}$ holds, except when $L = \FF_{q^2}$, in which case $L \cap C = F \cap C$.
\end{theorem}

Using Theorem \ref{MersenneException}, we prove the Main Theorem.

\begin{proof}[Proof of Main Theorem]
Due to It\^o's theorem (Theorem 6.15 of \cite{Isaacs}), every degree in $\cd(G)$ divides $\ind GC = |\varGamma| = e$.

Consider a character $\mu \in \Irr (C)$.  Let $\varPsi$ be the stabilizer of $\mu$ in $\varGamma$.  It follows that $C\varPsi$ is the stabilizer of $\mu$ in $G$.  As $C\varGamma/C \iso \varGamma$ is cyclic, we see that $\mu$ extends to $\Irr(C\varPsi)$ (see Corollary 11.22 of \cite{Isaacs}) and by Gallagher's theorem, every irreducible constituent of $\mu^{C\varPsi}$ is an extension of $\mu$.  Applying the Clifford correspondence (Theorem 6.11 of \cite{Isaacs}), we see that every irreducible constituent of $\mu^{C\varPsi}$ induces irreducibly to $G$.  We conclude that every irreducible constituent of $\mu^G$ has degree $|\varGamma:\varPsi|$.  Conversely, since every irreducible character of $G$ lies over some
character $\mu \in \Irr (C)$, every degree in $\cd (G)$ arises in this manner.


We now establish some notation to be used throughout the remainder of the proof.  Let $\lambda$ be a generator for the cyclic group $\Irr (C)$.  Note that $\lambda$ is a faithful character and a homomorphism from $C$ to the complex numbers.  Let $c$ be a generator for the group $C$, and $\sigma$ a generator for $\varGamma$.  (In particular, we may take $\sigma$ to be the Frobenius automorphism of $E$.)

Consider an arbitrary integer $m$ and element $\tau \in \varGamma$.  Because $\lambda$ is a homomorphism, we have
$$
\lambda^m (d) = \lambda (d^m)
$$
for all elements $d \in C$.  We compute
$$
(\lambda^m)^\tau (d) = \lambda^m (d^{\tau^{-1}}) = \lambda \big( (d^{\tau^{-1}})^m \big) = \lambda \big( (d^m)^{\tau^{-1}} \big).
$$
Thus, we have $(\lambda^m)^\tau = \lambda^m$ if and only if $\lambda \big( (d^m)^{\tau^{-1}} \big) = \lambda (d^m)$ for all $d \in C$.  Since $\lambda$ is faithful, this will occur if and only if $(d^m)^{\tau^{-1}} = d^m$ for all $d \in C$.  It follows that $\tau$ stabilizes $\lambda^m$ if and only if $\tau$ centralizes $d^m$ for every element $d \in C$.  The latter happens exactly when $\tau$ centralizes $\langle c^m \rangle$.

Suppose $n \ne 2$ is a divisor of $e$.  Let $\varPhi = \langle \sigma^n \rangle$, and observe that $\ind \varGamma\varPhi = n$.  Write $L$ for the fixed field of $\sigma^n$ in $E$.  Note that $n \not= 2$ implies that $L = \FF_{q^n} \ne \FF_{q^2}$, that $\varPhi = \Gal(E/L)$, and that $\sigma^n$ fixes the subgroup $L \cap C$ of~$C$.  Since $c$ is a generator of $C$, there is some integer $m$ for which $L \cap C = \langle c^m \rangle$.  We now compute the stabilizer of $\lambda^m$ in $\varGamma$.  Let $\tau \in \varGamma$ be arbitrary.  By the previous paragraph, we see that $\tau$ stabilizes $\lambda^m$ if and only if $\tau$ centralizes $\langle c^m \rangle = L \cap C$.  This occurs if and only if $\tau$ centralizes $\widehat {L \cap C}$.  Applying Theorem~\ref{MersenneException}, we know that $L = \widehat {L \cap C}$.  Thus, $\tau$ stabilizes $\lambda^m$ if and only if $\tau \in \Gal (E/L) = \varPhi$.  We conclude that $\varPhi$ is the stabilizer of $\lambda^m$, and applying the second paragraph, we have $n \in \cd (G)$.

To see that $2$ is not in $\cd (G)$, it suffices by the second paragraph to show that $\langle \sigma^2 \rangle$ is not the stabilizer of any character in $\Irr (C)$.  In other words, we need to prove that if $\sigma^2$ stabilizes $\lambda^l$ for some integer $l$, then $\sigma$ stabilizes $\lambda^l$.  In the fourth paragraph, we have seen that $\sigma^2$ stabilizes $\lambda^l$ if and only if $\sigma^2$ centralizes $c^l$.  Let $K$ be the fixed field for $\sigma^2$ in $E$.  It follows that $K$ is the extension field of $F$ whose degree is $2$.  (I.e., $K=\FF_{q^2}$.)  {}From Theorem~\ref{MersenneException}, we know that $K \cap C = F \cap C$.  Now, we see that $\sigma^2$ stabilizes $\lambda^l$ if and only if $c^l \in K \cap C = F \cap C$.  But $c^l \in F$ implies that $\sigma$ stabilizes $\lambda^l$, and the result is proved.
\end{proof}

\end{document}